\documentclass{article}

\usepackage{arxiv}

\usepackage[utf8]{inputenc} % allow utf-8 input
\usepackage[T1]{fontenc}    % use 8-bit T1 fonts
\usepackage{hyperref}       % hyperlinks
\usepackage{url}            % simple URL typesetting
\usepackage{booktabs}       % professional-quality tables
\usepackage{amsfonts}       % blackboard math symbols
\usepackage{nicefrac}       % compact symbols for 1/2, etc.
\usepackage{microtype}      % microtypography
\usepackage{lipsum}		% Can be removed after putting your text content
\usepackage{graphicx}
\usepackage{natbib}
\usepackage{doi}
\usepackage{amsmath, amsthm, amssymb}
\usepackage{mathrsfs}
\newcommand{\bra}[1]{\left(#1\right)}

\newcommand{\seq}[1]{\left<#1\right>}
\newcommand{\vertiii}[1]{{\left\vert\kern-0.25ex\left\vert\kern-0.25ex\left\vert #1
    \right\vert\kern-0.25ex\right\vert\kern-0.25ex\right\vert}}
\newcommand{\norm}[1]{\left\Vert#1\right\Vert}
\newcommand{\abs}[1]{\left\vert#1\right\vert}
\newcommand{\set}[1]{\left\{#1\right\}}
\renewcommand{\c}{\mathbb C}
\newcommand{\R}{\mathbb R}
\newcommand{\N}{\mathbb N}
\renewcommand{\r}{\mathrm{ran}}
\newcommand {\bh}{\mathscr{B}(\mathscr{H})}
\newcommand {\lh}{\mathscr{L}(\mathscr{H})}
\newcommand {\n}{\ker}
\newcommand {\D}{\mathfrak{D}}
\newcommand {\M}{\mathcal{M}}
\newcommand {\T}{\mathbb{T}}
\newcommand {\h}{\mathscr{H}}
\newcommand {\ind}{\mathrm{ind}}
\renewcommand {\L}{\mathscr{L}}
\renewcommand {\b}{\mathscr{B}}
\newtheorem{theorem}{Theorem}[section]
\newtheorem{lemma}[theorem]{Lemma}
\newtheorem{proposition}[theorem]{Proposition}
\newtheorem{corollary}[theorem]{Corollary}
\newtheorem{definition}[theorem]{Definition}
\newtheorem{example}[theorem]{Example}

\newtheorem{remark}[theorem]{Remark}
\newcommand\mystyle{\everymath{\displaystyle}}
\mystyle
\title{On the closability of class totally paranormal operators}

\author{Salam Alnabulsi$^{1}$, \href{https://orcid.org/0000-0002-3816-5287}{\includegraphics[scale=0.06]{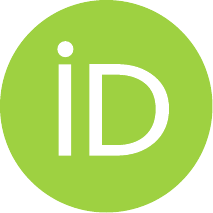}\hspace{1mm}M.H.M.~Rashid$^{2}$}\thanks{Corresponding Author} \\
$^{1}$ Department of Mathematics, Faculty of Science, University of Jordan, Amman 11942,
Jordan\\
	$^{2}$Department of Mathematics\&Statistics\\Faculty of Science P.O.Box(7)\\
	Mutah University University\\
	Mu'tah-Jordan \\
	\texttt{mrash@mutah.edu.jo}
}

\renewcommand{\shorttitle}{\textit{arXiv} Template}

\hypersetup{
pdftitle={On the closability of class totally paranormal operators},
pdfsubject={q-bio.NC, q-bio.QM},
pdfauthor={Salam Alnabulsi, M.H.M.Rashid},
pdfkeywords={Densely defined operator; closed operator; totally paranormal;
reduced minimum modulus; Riesz projection and Weyl’s theorem},
}

\begin{document}
\maketitle

\begin{abstract}
	This article delves into the analysis of various spectral properties pertaining to totally paranormal closed operators, extending beyond the confines of boundedness and encompassing operators defined in a Hilbert space. Within this class, closed symmetric operators are included.
Initially, we establish that the spectrum of such an operator is non-empty and provide a characterization of closed-range operators in terms of the spectrum. Building on these findings, we proceed to prove Weyl's theorem, demonstrating that for a densely defined closed totally paranormal operator $T$, the difference between the spectrum $\sigma(T)$ and the Weyl spectrum $\sigma_w(T)$ equals the set of all isolated eigenvalues with finite multiplicities, denoted by $\pi_{00}(T)$.
In the final section, we establish the self-adjointness of the Riesz projection $E_{\mu}$ corresponding to any non-zero isolated spectral value $\mu$ of $T$. Furthermore, we show that this Riesz projection satisfies the relationships $\mathrm{ran}(E_{\mu}) = \n(T-\mu I) = \n(T-\mu I)^*$. Additionally, we demonstrate that if $T$ is a closed totally paranormal operator with a Weyl spectrum $\sigma_w(T) = {0}$, then $T$ qualifies as a compact normal operator.
\end{abstract}

% keywords can be removed

\section{Introduction}
The class of normal operators stands out as a crucial and extensively researched category in operator theory. The spectral theorem for normal operators guarantees the presence of non-trivial invariant subspaces and unveils the comprehensive structure of the operator.
The class of bounded paranormal operators was first studied by Istr\v atescu \cite{Istr}, who
named it as the class $N$. Further, Furuta \cite{Furuta67} introduced the term paranormal operator.
Many authors studied bounded paranormal operators, for example \cite{Ando, Furuta67, Istr}.
In particular, Ando \cite{Ando} gave a characterization of bounded paranormal operators.
Istr\v atescu \cite{Istr} proved that the class of normaloid operators is a generalization of
paranormal operators.
A continuous linear operator on a complex Banach space is said to be paranormal if $\norm{Tx}^2\leq \norm{T^2x}\norm{x}$ for all $x\in X$, where $X$ is a Banach space. $T$ is called totally paranormal \cite{Schmo} if $T-\mu I$ is paranormal
for every $\mu\in\c$. That is, $\norm{(T-\mu I)x}^2\leq \norm{(T-\mu I)^2x}\norm{x}$ for all $x\in  X$ and $\mu\in\c$. Hence,
We have the following inclusion relation between some subclasses and a generalized class of bounded totally paranormal operators.
\begin{equation*}
  \mbox{Normal}\subseteq \mbox{Hyponormal}\subseteq \mbox{Totally Paranormal}\subseteq \mbox{Paranormal}\subseteq \mbox{normaloid}.
\end{equation*}
The above inclusion relations are proper. For more details, we refer to \cite{Furuta67, Schmo}. The
definition of bounded paranormal operators is extended to unbounded operators by
Daniluk \cite{Daniluk}, where he discussed about closability of unbounded paranormal operators.\\
\indent In this article, we are going to deal with densely defined closed totally paranormal operators
in a Hilbert space $\h$ and prove the following results.\\
\indent Let $T$ be a densely defined closed totally  paranormal operator in $\h$. Then
\begin{itemize}
  \item [(i)] spectrum of $T$ is non-empty.
  \item [(ii)] Every isolated spectral value of $T$ is an eigenvalue.
  \item [(iii)] In addition, if $\n(T)= \n(T^*)$, then
  \begin{enumerate}
   \item[(a)] range of $T$ is closed if and only if 0 is an isolated spectral value of $T$.
    \item[(b)] The minimum modulus, $m(T)$ is equal to the distance of 0 from spectrum of $T$.
  \end{enumerate}
  \item [(iv)] $T$ satisfies the Weyl's Theorem i.e. $\sigma(T)\setminus\sigma_w(T)=\pi_{00}(T)$. Here $\sigma_w(T)$ is the
Weyl's spectrum and $\pi_{00}(T)$ consists of all isolated eigenvalues of $T$ with finite
multiplicity.
\item [(v)] If $\mu$ is a non-zero isolated spectral value of $T$, then the Riesz projection $E_{\mu}$ with
respect to $\mu$ is self-adjoint and satisfies $\r(E_{\mu})=\n(T-\mu I)=\n(T-\mu I)^*$.
\end{itemize}
\indent The exploration of Weyl's theorem and the self-adjointness of the Riesz projection concerning an isolated spectral value has been undertaken across various classes of operators. Coburn \cite{Coburn} established these properties for certain non-normal operators, specifically hyponormal and Toeplitz operators. Building upon this, Schmoeger \cite{Schmo} extended the investigation to bounded totally paranormal operators, utilizing Ando's characterization \cite{Ando} for paranormal operators. However, Ando's characterization is not applicable to unbounded paranormal operators, and the techniques employed for bounded operators are not suitable in this context. Consequently, we endeavor to prove properties (iv) and (v) through an alternative approach.\\
\indent Gupta and Mamtani \cite{Gup-Mam} demonstrated the fulfillment of Weyl's theorem by closed hyponormal operators. In their subsequent work \cite{Gup-Mam-1}, the same authors presented several essential conditions that must be met for the orthogonal direct sum of densely defined closed operators to satisfy Weyl's theorem.\\
\indent It seems like you've provided a statement about an article, possibly related to functional analysis or operator theory. If you have a specific question or if there's a particular aspect of the article you'd like to discuss or elaborate on, please provide more details or ask a specific question. I'll do my best to assist you based on the information you provide.\\
\indent The article is structured into four sections for clarity. In the second section, we establish essential notations and highlight pertinent known results that will be employed consistently throughout the paper. Moving on to the third section, we delve into an exploration of various spectral properties inherent in densely defined closed totally paranormal operators. Finally, the fourth section is dedicated to the proof of Weyl's theorem specifically for densely defined closed totally paranormal operators.

%=================================================================================
\section{Notations and preliminaries}
%==============================================================================
In this article, we explore intricate Hilbert spaces, represented as $\h, \h_1, \h_2$, and so forth. The inner product and the corresponding norm are symbolized by $\seq{\cdot,\cdot}$ and $\norm{\cdot}$, respectively.\\
\indent The set of all linear operators on $\h$ is denoted as $\lh$, while the collection of all bounded linear operators is represented as $\bh$. For a linear operator $T \in \lh$, we use $\D(T)$, $\n(T)$, and $\r(T)$ to signify its domain, null space, and range space, respectively. A linear operator $T$ is termed a \emph{densely defined operator} if $\overline{\D(T)}=\h$.\\
\indent If $T\in\lh$ and $\M$ is a closed subspace of $\h$ , then $\M$ is said to be invariant under
$T$ , if for every $x\in\D(T)\cap \M$, $Tx$ is in $\M$. We denote the identity operator on $\M$ by
$I_{\M}$, the orthogonal projection on $\M$ by $P_{\M}$. The unit sphere of $\M$ is $\T_{\M}:=\set{x\in\M: \norm{x}=1}$.
 The restriction of $T$ to $\M$ is an operator $T|_{\M}:\M\cap\D(T)\to \h$
 defined  by $T|_{\M}x=Tx$,
 for all $x\in \M\cap D(T)$.
 If $\M$ is invariant under $T$ , then $T|_{\M}$
is an operator from $\D(T)\cap \M$ into $\M$.\\
\indent  An operator $T\in\lh$ is said to be \emph{closed} if for any sequence $\{x_n\}\subseteq \D(T)$ with $x_n\to x$
 and $Tx_n\to y$ then $x\in\D(T)$ and $Tx=y$. In this document, the notation $C(\h_1,\h_2)$ will be employed to denote the collection of closed linear operators such that $\D(T)\subseteq \h_1$ and $\r(T)\subseteq \h_2$. In the case where $\h_1$ equals $\h_2$, we will use the shorthand $C(\h)$.\\
 It is known that every densely defined operator $T\in C(\h_1,\h_2)$ has a unique adjoint in
  $ C(\h_2,\h_1)$, that is, there exists a unique $T^*\in C(\h_2,\h_1)$  such that $\seq{Tx,y}=\seq{x,T^*y}$ for all $x\in\D(T)$ and $y\in\D(T^*)$.
%===========================================================================================
 \begin{remark} By the closed graph theorem (cf. \cite[Theorem 7.1, page 231]{Nair}, it follows that
a closed operator $T\in C(\h_1,\h_2)$ with $\D(T)=\h_1$ is bounded.
 \end{remark}
%==============================================================================
 \begin{lemma}\cite{Ku-Ra}\label{Lemma1} Let $T\in\lh$ be a densely  defined closed operator. Then
 \begin{equation*}
   \overline{\D(T)\cap\n(T)^{\bot}}=\n(T)^{\bot}.
 \end{equation*}
 \end{lemma}
 %==============================================================================
 If $S$ and $T$ are two closed operators, then $S$ is called an extension of $T$ (or $T$ is a
restriction of $S$), if $\D(T )\subseteq \D(S)$ and $Sx=Tx$ for all $x\in\D(T)$. This is often
denoted as $T\subseteq S$. Consequently, $S=T$ if and only if $\D(S)=\D(T)$ and $Sx=Tx$ for all
$x\in \D(S)=\D(T)$.
%========================================================================================
\begin{definition}\cite{Walter} A densely defined operator $T\in\lh$ is said to be self-adjoint if $\D(T)=\D(T^*)$
and $T=T^*$. And a self-adjoint $T$ is said to be positive if $\seq{Tx,x}\geq 0$ for all $x\in\D(T)$.
\end{definition}
%========================================================================================
\begin{definition}\cite[Page 365]{Rudin} If $T\in\lh$ is a closed operator, then the resolvent set
of $T$ is defined by
\begin{equation*}
  \rho(T)=\set{\mu\in\c:T-\mu I\,\,\mbox{is invertible and $(T-\mu I)^{-1}\in\bh$}}
\end{equation*}
and the spectrum of $T$, denoted by $\sigma(T)$, is defined by
\begin{equation*}
  \sigma(T):=\c\setminus\rho(T)
\end{equation*}
\end{definition}
Note that $\sigma(T)$ is a closed subset of $\c$. Moreover $\sigma(T)$ can be empty set or the whole
complex plane $\c$.\\
The spectrum of $T$ decomposes as the disjoint union of \emph{the point spectrum} $\sigma_p(T)$,
\emph{the continuous spectrum} $\sigma_c(T)$ and \emph{the residual spectrum} $\sigma_r(T)$, where
\begin{eqnarray*}
  \sigma_p(T) &=&\set{\mu\in\c:T-\mu I\,\,\mbox{is not injective}}, \\
  \sigma_r(T) &=&\set{\mu\in\c:T-\mu I\,\,\mbox{is not injective but $\r(T-\mu I)$ is not dense in $\h$}},\\
  \sigma_c(T)&=&\sigma(T)\setminus\bra{\sigma_p(T)\cup\sigma_r(T)}.
\end{eqnarray*}
The spectral radius of $T\in\bh$ is defined by
\begin{equation*}
  r(T):=\sup\set{|\mu|: \mu\in\sigma(T)}.
\end{equation*}
An operator $T\in\bh$ is said to be normaloid, if $r(T)=\norm{T}$.\\
Recall that a linear operator $T\in\lh$ is compact, if $T$ maps every bounded set in $\h$
to a pre-compact set in $\h$. For more details about compact  operators,
we refer to \cite{Schechter}.
%====================================================================
\begin{definition}\cite[Page 156]{Schechter} A closed operator $T$ in a densely defined space $\h$ is termed Fredholm if the range $\r(T)$ is closed, and both the dimensions of the null space $\n(T)$ and its orthogonal complement $\r(T)^\perp$ are finite. In such instances, the index of $T$, denoted by $\mathrm{ind}(T)$, is defined as $\ind(T)=\dim(\n(T))-\dim(\r(T)^{\bot})$.
\end{definition}
\begin{remark}\label{remark1-intro} If $T\in\lh$ is a densely defined closed Fredholm operator and $K$ is a
compact operator, then $T+K$ is also Fredholm and $\ind(T+K)=\ind(T)$.
\end{remark}
%=======================================================================
\begin{definition}\cite[Page 172]{Schechter} If $T\in\lh$ is a densely defined closed operator, then
the Weyl's spectrum of $T$ is defined by
\begin{equation*}
  \sigma_w(T)=\set{\lambda\in\c:T-\lambda I\,\mbox{is not Fredholm of index 0}}
\end{equation*}
and  $\pi_{00}(T)=\set{\lambda\in\sigma_p(T):\lambda\,\,\mbox{is isolated with $\dim\bra{\n(T-\lambda I)}<\infty$}}$.
\end{definition}
Suppose $T\in\lh$ is a densely defined closed operator with $\sigma(T)=\sigma\cup \tau$, where $\sigma$
is contained in some bounded domain $\Delta$ such that 	$\overline{\Delta}\cap \tau=\emptyset$. Let $\Lambda$
 be the boundary of $\Delta$, then
 \begin{equation}\label{Equation2.1}
   E_{\sigma}=\frac{1}{2\pi i}\int_{\Lambda}(z I-T)^{-1}dz
 \end{equation}
 is called the Riesz projection with respect to $\sigma$.
%=======================================================================================
\begin{theorem}\cite[Theorem 2.1, Page 326]{GGK} \label{Theorem1}Suppose $T\in\lh$ is a densely defined
closed operator with $\sigma(T)=\sigma\cup \tau$, where $\sigma$ is contained in some bounded domain
	 and $E_{\sigma}$ is the operator defined in Equation (\ref{Equation2.1}). Then
\begin{enumerate}
  \item [(i)] $E_{\sigma}$ is a projection.
  \item [(ii)] The subspaces $\r(E_{\sigma})$ and $\n(E_{\sigma})$ are invariant under $T$.
  \item [(iii)] The subspace $\r(E_{\sigma})$ is contained in $\D(T)$ and $T|_{\r(E_{\sigma})}$ is bounded.
  \item [(iv)] $\sigma\bra{T|_{\r(E_{\sigma})}}=\sigma$ and $\sigma(T|_{\n(E_{\sigma})})=\tau$.
\end{enumerate}
\end{theorem}
In particular, if $\mu$ is an isolated point of $\sigma(T)$, then there exist a positive real number
$r$ such that $\{z\in\c: |z-\mu|\leq r\}\cap \sigma(T)=\set{\mu}$. If we take $\set{z\in\c:|z-\mu|=r}$,
then the Riesz projection with respect to $\mu$ is defined by
\begin{equation}\label{Eq.2.2}
  E_{\mu}=\frac{1}{2\pi i}\int_{\Lambda}(z I-T)^{-1}dz.
\end{equation}
%====================================================================================
\begin{definition}\cite{Ku-Ra-1} Let $T\in\lh$ be a closed operator. Then
\begin{enumerate}
  \item [(i)] the minimum modulus of $T$ is defined by $m(T):=\set{\norm{Tx}: x\in \T_{\D(T)}}$. Then
  \item [(ii)]  the reduced minimum modulus of $T$ is defined by $\gamma(T):=\set{\norm{Tx}:x\in\T_{\D(T)\cap\n(T)^{\bot}}}$.
\end{enumerate}
By the definition, it is clear that  $m(T)\leq \gamma(T)$
\end{definition}
%============================================================================
The following characterization of closed range operators is frequently used in the
article.
\begin{theorem}\cite[Page 334]{Ben-Israel} \label{Theorem2} For a densely defined closed operator $T\in\lh$, the
following are equivalent.
\begin{enumerate}
  \item [(i)] $\r(T)$ is closed.
  \item [(ii)] $\r(T^*)$ is closed.
  \item [(iii)] $\gamma(T)>0$.
  \item [(iv)] $S_0=T|_{\D(T)\cap\n(T)^{\bot}}$ has a bounded inverse.
\end{enumerate}
\end{theorem}
If $T\in\lh$ is a densely defined closed operator and  $\n(T)=\set{0}$, then the
inverse operator, $T^{-1}$ is the linear operator from $\h$ to $\h$, with $\D(T^{-1})=\r(T)$ and
$T^{-1}Tx=x$ for all $x\in \D(T)$. In particular if $T$ is a bijection, then by the closed
graph theorem it follows that $T^{-1}\in\bh$. In addition, if $T$ is normal then $T$ has a
bounded inverse if and only if $m(T)> 0$.
%===============================================================
\begin{theorem}\label{Moh-closed}\label{Theorem3} If $T\in\bh$ is a totally paranormal, then
\begin{enumerate}
  \item [(i)] $T$ is normaloid.
  \item [(ii)] $T^{-1}$ is totally paranormal, if $T$ is invertible.
  \item [(iii)] $T$ is unitary, if $\sigma(T)$ lies on the unit circle.
\end{enumerate}
\end{theorem}
%===============================================================
%=======================================================================================
\section{Spectral Properties}
%====================================================
In this section, we study some spectral properties of densely defined closed totally paranormal
operators.
\begin{definition} A densely defined operator $T\in\lh$ is said to be totally paranormal
if
\begin{equation*}
  \norm{(T-\lambda I)x}^2\leq \norm{(T-\lambda I)^2x}\norm{x}\,\,\mbox{ for all $\lambda\in\c$ and} \,\,\D((T-\lambda)^2)\subseteq \D(T-\lambda).
\end{equation*}
Equivalently, $T$ is totally paranormal if and only if $T-\lambda I$ is paranormal for all $\lambda\in\c$. And $T$ is totally $*$-paranormal
if $T-\lambda I$ is $*$-paranormal for all $\lambda\in\c$.
\end{definition}
A densely defined operator $T$ in $\h$ is said to be hyponormal if $\D(T)\subseteq \D(T^*)$
and $\norm{T^*x}\leq \norm{Tx}$ for $x\in\D(T)$. And $T$ is cohyponormal if $T^*$ is hyponormal.
%====================================================================
\begin{proposition} Let $T\in C(\h)$. Then
\begin{enumerate}
  \item [(i)] If $T$ is closed cohyponormal and closed totally $*$-paranormal, then
  $T$ is closed totally paranormal.
  \item [(ii)] If $T$ is closed hyponormal and closed totally paranormal, then
  $T$ is closed totally $*$-paranormal.
\end{enumerate}
\end{proposition}
\begin{proof} (i) Assume that $T$ is closed cohyponormal and closed totally $*$-paranormal.
 If $T$ is closed cohyponormal, then so is $T-\lambda I$ for all $\lambda\in\c$. Hence,
 \begin{equation*}
   \D((T-\lambda I)^2)\subseteq \D((T-\lambda I)^*)\subseteq \D(T-\lambda I).
 \end{equation*}
 Let $x\in\D((T-\lambda I)^2)$. Then we have
 \begin{equation*}
   \norm{(T-\lambda I)x}^2\leq \norm{(T-\lambda I)^*x}^2\leq \norm{(T-\lambda I)^2x}\norm{x}.
 \end{equation*}
 (ii) By the hypotheses on $T$, $\D((T-\lambda I)^2)\subseteq \D(T-\lambda I)\subseteq \D((T-\lambda I)^*)$, and for each
 $x\in \D((T-\lambda I)^2)$,
 \begin{equation*}
   \norm{(T-\lambda I)^*x}^2\leq \norm{(T-\lambda I)x}^2\leq \norm{(T-\lambda I)^2x}\norm{x}.
 \end{equation*}
 So, the proof is complete.
\end{proof}
%=====================================================================
\begin{proposition} Let $T\in C(\h)$ be totally $*$-paranormal, and let $\lambda$ be any complex scalar.
Then $\n(T-\lambda I)\subset \n(T-\lambda I)^* $.
\end{proposition}
\begin{proof} Let $x\in \D(T)$ be a unit eigenvector of $T$  associated to $\lambda$. Then $Tx=\lambda x$.
By the hypotheses on $T$, $\norm{T^*x}\leq \abs{\lambda}$. Consequently,
\begin{eqnarray*}
  0 &\leq&\norm{T^*x-\bar{\lambda}x}^2=\norm{T^*x}^2-\lambda\seq{x,Tx}-\bar{\lambda}\seq{Tx,x}+|\lambda|^2 \\
   &=&\norm{T^*x}^2-|\lambda|^2\leq 0.
\end{eqnarray*}
Thus, $T^*x=\bar{\lambda}x$ and so $x\in\n(T-\lambda)^*$.
\end{proof}
%==========================================================================
\begin{lemma}\label{lemma-Mohad} Let $T\in C(\h)$ be a totally $*$-paranormal, then  there exists a contraction $Q_{\lambda}\in\bh$
    such that $(T-\lambda I)^2\subset (T-\lambda I)^* Q_{\lambda}$.
\end{lemma}
\begin{proof} Since  $\norm{(T-\lambda I)^*x}^2\leq \norm{(T-\lambda I)^2x}\norm{x}$ for all $x\in\D((T-\lambda I)^2)\subseteq \D((T-\lambda
I)^*)$, there exists a contraction $K'\in\b\bra{\overline{\r\bra{(T-\lambda I)^2}},\overline{\r(T^*-\bar{\lambda}I)}}$ such that
$K'(T-\lambda I)^2\subset T^*-\bar{\lambda}I$. Let $K\in\bh$ be any contraction which extends $K'$ (e.g. set $Kx=0$ for
 $x\in\h\ominus \overline{\r\bra{(T-\lambda I)^2}}$). Then $K(T-\lambda I)^2\subset T^*-\bar{\lambda}I$.  Taking
adjoints in the last inclusion and exploiting the closability of $T$ , we get $(T-\lambda I)^2\subseteq \bra{(T-\lambda I)^2}^{**}
\subseteq \bra{K(T-\lambda I)^2}^*=(T-\lambda I)^{*2}K^*\subseteq (T-\lambda I)^*K^*$. This gives us the conclusion with  $Q=K^*$.
\end{proof}
%============================================================================
A closed subspace $\M$ of $\h$ reduces $T\in C(\h)$ if $\M$ and $\M^{\bot}$ are invariant
under $T$. Stochel \cite{Stochel} proved if $T\in C(\h)$ is hyponormal and $\M$ is a closed
subspace of $\h$ which is invariant under $T$ with $T|_{\M}$ is normal, then $\M$ reduces
$T$.
\begin{theorem} Suppose $T\in C(\h)$ is a densely defined totally $*$-paranormal. If $\M$ is a closed subspace of $\h$
which is invariant under $T$ with $T|_{\M}$ is
normal, then $\M$ reduces $T$.
\end{theorem}
\begin{proof} Let $\h=\h_1\oplus \h_2$, where $\h_1=\M$ and $\h_2=\M^{\bot}$. Then $T$ has the
block matrix representation
\begin{equation*}
  \begin{bmatrix}
  T_{11} & T_{12} \\
  T_{21} & T_{22} \\
\end{bmatrix},
\end{equation*}
where $T_{ij}: \D(T)\cap \h_{j}\to \h_i$ is defined by $T_{ij}=P|_{\h_{i}}TP|_{\h_j}|_{\D(T)\cap H_{k}}$ for $k=1,2$.
Here, $P|_{\h_i}$ denotes the orthogonal projection onto $\h_i$. Since $\M$ is invariant
under $T$, we have
  \begin{equation*}
  \begin{bmatrix}
  T_{11} & T_{12} \\
  0 & T_{22} \\
\end{bmatrix}.
\end{equation*}
Let $y\in \D(T)\cap\M^{\bot}$. By Lemma \ref{lemma-Mohad}, we have
\begin{equation*}
  (T-\lambda I)^2\subset (T-\lambda I)^* Q_{\lambda}
\end{equation*}
for every $\lambda\in\c$. Thus, $\r\bra{(T-\lambda I)^2}\subseteq \r\bra{(T-\lambda I)^*}$for every $\lambda\in \c$. Then
there exist a densely defined operator $B$ such that $(T-\lambda I)^2=(T-\lambda I)^*B$ (see \cite{Douglas}). Hence, $T_{12}(y) = (T_{11}-\lambda I)^*u$
 for some $u\in \M$. We can choose $v$ such that $(T_{11}-\lambda I)^*u=(T_{11}-\lambda I)v$.
  Therefore, $T_{12}(y)=(T_{11}-\lambda I)v$ for every $\lambda\in\c$. Consequently,
  \begin{equation*}
    T_{12}(y)\in\bigcap_{\lambda\in\c}\r(T_{11}-\lambda I).
  \end{equation*}
Hence, $T_{12}(y)=0$ for all $y\in \D(T)\cap\M^{\bot}$ (see \cite{Putnam}) and so $T_{12}=0$. This ends the proof.
\end{proof}
%====================================================================
The ascent $p(T)$ and descent $q(T)$ of an operator $C(\h)$ are given by
\begin{eqnarray*}
  p(T) &=& \inf\set{n :\n(T^n)=\n(T^{n+1})}\,\,\mbox{and} \\
  q(T) &=&\inf\set{n:\r(T^n)=\r(T^{n+1})}
\end{eqnarray*}
It follows from the definition of totally paranormal the following result holds.
\begin{proposition} Let $T\in\lh$ be a densely defined closed totally paranormal operator.
Then the ascent $p(T)$ and descent $q(T)$ of $T$ are finite for all $\lambda\in\c$.
\end{proposition}
%%%%%%%%%%%%%%%%%%%%%%%%%%%%%%%%%%%%%%%%%%%%%%%%%%%%%%%%%%%%%%%%%%%%%%%%%%%%%%%%%%%%%%%
\begin{definition}\cite{Kato} Let $T$ be a non necessarily bounded operator with domain $\D(T) \subset \h$.
We say that $\lambda$ is not in $\sigma(T)$ if $T-\lambda$ is injective and $(T-\lambda)^{-1}\in\bh$.
\end{definition}
%=================================================================
The following results  immediately follows from the definition.
%======================================================================================
\begin{proposition} Let $T\in\lh$ be a densely defined closed totally paranormal operator. Then $T-\alpha I$ and
$\alpha T$ are totally paranormal operators for all $\alpha\in\c$.
\end{proposition}
%============================================================================
\begin{proposition} Let Let $T\in\lh$ be a densely defined closed totally paranormal operator.
If $\sigma(T)=\{\mu\}$, then $T=\mu I$.
\end{proposition}
%=================================================================================
We now give the  counterexample of a closed densely
defined operator $T$ such that both $T$ and $T^*$ are one-to-one and  totally paranormal, yet $T$  is not normal.
\begin{example} The Hilbert space in question is $L^2(\R)\oplus L^2(\R)$. From \cite{ota},
we have an explicit example of a densely defined unbounded closed
operator $T$ for which: $\D(T^2)=\D(T^{*2})=\{0\}$
More precisely, $T$  is defined by
\begin{equation*}
  T=\begin{bmatrix}
      0 & A^{-1} \\
      B & 0 \\
    \end{bmatrix}
\end{equation*}
  on $\D(T):=\D(B)\oplus \D(A^{-1})\subset L^2(\R)\oplus L^2(\R)$, and where $A$ and $B$ are
two unbounded self-adjoint operators such that $\D(A)\cap \D(B)=\D(A^{-1})\cap \D(B^{-1})=\{0\}$,
where $A^{-1}$ and $B^{-1}$ are not bounded (as in \cite{Kosaki}). Hence
\begin{equation*}
  T^*=\begin{bmatrix}
      0 & B \\
      A^{-1} & 0 \\
    \end{bmatrix}
\end{equation*}
for $A^{-1}$ and $B$ are both self-adjoint. Observe now that both $T$ and $T^*$
are one-to-one since both $A^{-1}$ and $B$ are so.
Both $T$  and $T^*$  are trivially totally paranormal thanks to the assumption
$\D(T^2)=\D(T^{*2})=\{0\}$. So tatally paranormality of both operators need only
be checked at the zero vector and this is plain as
$\norm{(T-\mu I)x}^2=\norm{(T-\mu I)^2x}\norm{x}=0$ and  $\norm{(T-\mu I)^{*}x}^2=\norm{(T-\mu I)^{2}x}\norm{x}=0$
for $x=0$. However, $T$ cannot be normal for it were, $T^2$   would be normal
too, in particular it would be densely defined which is impossible
here
\end{example}
%=====================================================================
Here we discuss some basic results related to unbounded totally paranormal operators, which
are often used in the article.

\begin{theorem} \label{theorem4}Let $T\in\lh$ be a densely defined closed totally paranormal operator and not a multiple of the identity. Then
the following holds.
\begin{enumerate}
  \item [(i)] If $\M$ is a closed invariant subspace of $T$ , then $T|_{\M}$ is totally paranormal.
  \item [(ii)] If $0\notin \sigma(T)$, then $T^{-1}$ totally paranormal
  \item [(iii)] $\sigma(T)$ is nonempty.
\end{enumerate}
\end{theorem}
\begin{proof}(i) For every $\lambda\in\c$, as $\M$ is invariant under $T$, we have
\begin{eqnarray*}
  \D((T-\lambda I)^2|_{\M}) &=&\D\bra{(T-\lambda I)^2}\cap \M \\
   &=& \set{x\in\D(T-\lambda I):(T-\lambda I)x\in \D(T-\lambda I)}\cap \M\\
   &=&\set{x\in\D(T-\lambda I)\cap \M:(T-\lambda I)x\in \D(T-\lambda I)\cap\M}\\
   &&\bra{\mbox{since $T(\D(T-\lambda I)\cap \M)\subseteq \M$}}\\
   &=&\set{x\in\D((T-\lambda I)|_{\M}):Tx\in \D((T-\lambda I)|_{\M}) }\\
   &=&\D\bra{\bra{(T-\lambda I)|_{\M}}^2}.
\end{eqnarray*}
Thus, $(T-\lambda I)^2|_{\M}=\bra{(T-\lambda)|_{\M}}^2$. Now the result follows from the below inequality;
\begin{eqnarray*}
  \norm{(T-\lambda I)|_{\M}x}^2 &=& \norm{(T-\lambda I)x}^2\leq \norm{(T-\lambda I)^2x} \\
   &=& \norm{(T-\lambda I)^2|_{\M}x}=\norm{\bra{(T-\lambda I)|_{\M}}^2x}, \forall x\in\T_{\D\bra{\bra{(T-\lambda I)|_{\M}}^2}}.
\end{eqnarray*}
(ii) Existence of $T^{-1}$ implies $\r(T)=\h$ and consequently $\r\bra{(T)^2}=\h$.
As $T$ is totally paranormal, we get $\n(T)=\n((T)^2)$, so $T^2$ is bijective and $\bra{(T)^2}^{-1}$ exists.
Also $\D(\bra{(T)^{-1}}^2)=\h=\r\bra{(T)^2}$. If $y \in \h$, then there exist $x\in \D\bra{(T)^2}$, such that
$y=T^2x$. Now,
\begin{eqnarray*}
  \norm{T^{-1}y}^2 &=& \norm{Tx}^2\leq \norm{T^2x}\norm{x} \\
   &=& \norm{y}\norm{T^{-2}y}.
\end{eqnarray*}
Hence $T^{-1}$ is totally paranormal since totally paranormal has invariant translation property.\\
(iii) Suppose  on the contrary that $\sigma(T)=\emptyset$. Then $T$ is invertible and
$T\in\bh$.\\
First, we show that $\sigma(T^{-1})=\{0\}$. For any complex number $\mu\neq 0$, consider
the operator $S=\mu^{-1} (T-\mu^{-1} I)^{-1}$. Here $S$ can also be written as the sum of
two bounded operators, $S=\mu^{-1}(I + \mu^{-1}(T-\mu^{-1} I)^{-1})$, so $S$ is bounded. By a
simple computation we can show that $S$ is the bounded inverse of $\mu I-T^{-1}$. Thus
$\sigma(T^{-1})\subseteq \{0\}$. As $T^{-1}\in\bh$, this implies $\sigma(T^{-1})$ is non-empty, so we conclude
that $\sigma(T^{-1})=\{0\}$.
By (ii), $T^{-1}$ is bounded totally paranormal operator and consequently normaloid by
Theorem \ref{Moh-closed}. Hence $\norm{T^{-1}}=0$, which implies $T^{-1}=0$, a contradiction. Hence
$\sigma(T)$ is non-empty.
\end{proof}
%===============================================================
Now we discuss about isolated spectral values of totally paranormal operators.
%=============================================================================
\begin{theorem}\label{Theorem5} Let $T$ be a densely defined closed totally paranormal operator. If $\mu$ is an
isolated point of $\sigma(T)$, then $\n(T-\mu I)=\r(E_{\mu})$.
\end{theorem}
\begin{proof} It follows from \cite[Lemma 3.4]{BR} that $\n(T-\mu I)\subseteq \r(E_{\mu})$. To complete the proof we have to show
that $\n(T-\mu I)\supseteq\r(E_{\mu})$. \\
\indent As a consequence of Theorem \ref{Theorem1} and Theorem \ref{theorem4}, we know that $T|_{\r\bra{E_{\mu}}}$ is
bounded and totally paranormal. By Theorem \ref{Theorem3} it follows that $T|_{\r\bra{E_{\mu}}}$ is normaloid.\\
If $\mu=0$, then     $\sigma\bra{T|_{\r(E_0)}}=\{0\}$. This implies $\norm{T|_{\r(E_0)}}=0$ and consequently
$T|_{\r(E_0)}=0$. Hence $\r(E_0)\subseteq \n(T)$.\\
If $\mu\neq 0$, then $\sigma\bra{\mu^{-1}T|_{\r\bra{E_{\mu}}}}=\{1\}$. By Theorem \ref{Theorem3}, it follows that $\mu^{-1}T|_{\r(E_{\mu})} $
is unitary. Thus $T|_{\r(E_{\mu})}-\mu I_{\r(E_{\mu})}$ is normal and $\sigma\bra{T|_{\r(E_{\mu})}-\mu I_{\r(E_{\mu})}}=\{0\}$. Since
every normal operator is normaloid, we conclude that $T|_{\r(E_{\mu})}-\mu I_{\r(E_{\mu})}=0$. Hence
$\r(E_{\mu})\subseteq \n\bra{T-\mu I}$.
\end{proof}\
%================================================================================
\begin{theorem}\label{Theorem6} Let $T\in\lh$ be a densely defined closed totally paranormal operator and
$\mu$ be an isolated point of $\sigma(T)$. Then $\n(E_{\mu})=\r(T-\mu I)$.
\end{theorem}
\begin{proof} By Theorem  \ref{Theorem1}, $\mu\notin \sigma\bra{T|_{\n(E_{\mu})}}$. This implies that $\r(T-\mu I)|_{\n(E_{\mu})}=n(E_{\mu})$
 and consequently $n(E_{\mu})\subseteq \r(T-\mu I)$.\\
 Let $y\in\r(T-\mu I)$. There exist $x\in \D(T)$ such that $y=(T-\mu I)x$. Since
 $\h=\r(E_{\mu})+\n(E_{\mu})$ and $\r(E_{\mu})\cap\n(E_{\mu})=\{0\}$, we have $x=p+q$, where $p\in \r(E_{\mu})$ and $q\in\n(E_{\mu})$.\\
 It follows from Theorem \ref{Theorem5}, that $p \in  \n(T-\mu I)\subseteq \D(T)$ and consequently
$q=x-p \in \D(T)$. As we know from Theorem \ref{Theorem1} that $\n(E_{\mu})$ is invariant under $T$,
we have
\begin{equation*}
  y=(T-\mu I)x=(T-\mu I)q \in (T-\mu I)(\n(E_{\mu}))\subseteq \n(E_{\mu}).
\end{equation*}
Hence $\r(T-\mu I)\subseteq \n(E_{\mu})$. This proves the result.
\end{proof}
%===============================================================================
The following results are consequences of Theorem \ref{Theorem6} which gives a characterization for closed range totally paranormal operators.
\begin{corollary}\label{Corollary1} Suppose $T\in\lh$ is a densely defined closed totally paranormal operator. If
0 is an isolated point of $\sigma(T)$, then $\r(T)$ is closed.
\end{corollary}
%=========================================================================
In general the converse of Corollary \ref{Corollary1} is not true. We have the following example to
illustrate this.
\begin{example}\label{example1} Let $T:\ell^2(\N)\to\ell^2(\N)$ be defined by
\begin{equation*}
  T(x_1,x_2,\cdots)=(0,x_1,x_2,\cdots),\,\,\mbox{for all $(x_n)\in\ell^2(\N)$}
\end{equation*}
Then $\sigma(T)=\{z\in\c:|z|\leq 1\}$, $\r(T)=\ell^2(\N)\setminus span\set{e_1}$. Here $\r(T)$ is closed but
0 is not an isolated point of $\sigma(T)$. Clearly, $T$ is a totally paranormal operator.
\end{example}
%=========================================================================================
Next result gives a sufficient condition under which the converse of Corollary \ref{Corollary1} is
also true.
\begin{theorem}\label{Theorem7} Let $T\in\lh$ be a densely defined closed totally paranormal operator with
$\n(T)=\n(T^*)$ and $0\in\sigma(T)$. Then $0$ is an isolated point of $\sigma(T)$ if and only if $\r(T)$
is closed.
\end{theorem}
\begin{proof} The necessary condition follows from Corollary \ref{Corollary1}. To prove the sufficient condition. Assume that
$\r(T)$ is closed. Consider $S_0=T|_{\n(T)^{\bot}}: \n(T)^{\bot}\cap\D(T)\to\n(T)^{\bot}$.
Clearly $S_0$ is injective and $\r(S_0)=\r(T )$ is closed. Also $\r(S_0)=\n(T^*)^{\bot}=\n(T)^{\bot}$,
 consequently $S_0$ is bijective and $S_0^{-1}\in \b(\n(T)^{\bot})$.
Thus $0\notin\sigma(S_0)$. Applying \cite[Theorem 5.4, Page 289]{Taylor}, $\sigma(T)\subseteq\{0\}\cup \sigma(S_0)$. Since
$0\in\sigma(T)$, we have $\sigma(T)=\{0\} \cup \sigma(S_0)$ and hence 0 is an isolated point of $\sigma(T)$.
\end{proof}
%=====================================================================================
Note that Theorem \ref{Theorem7} does not hold if we drop the condition $\n(T)=\n(T^*)$.
Consider the operator $T$ defined in Example \ref{example1}. Clearly $\n(T)=\{0\}\neq span\{e_1\}=\n(T^*)$,
and $\r(T)$ is closed but 0 is not an isolated point of $\sigma(T)$.
%=================================================================================
\begin{theorem}\label{theorem8} Let $T\in\lh$ be a densely defined closed totally paranormal operator. If
$\n(T)=\n(T^*)$, then $m(T)=d(0,\sigma(T))$, the distance between 0 and $\sigma(T)$.
\end{theorem}
\begin{proof} We will prove this result by considering the following two cases, which exhaust
all the possibilities.\\
Case (1): $T$ is not injective. Clearly $m(T)=0$ and $0\in\sigma_p(T)$. Hence $m(T)=0=d(0,\sigma(T))$.\\
Case (2): $T$ is injective. It suffices to show that $\gamma(T)=d(0,\sigma(T))$ because $m(T )=\gamma(T)$.\\
\indent First assume that $\gamma(T)=0$. It follows from Theorem \ref{Theorem2} that $\r(T)$ is not closed
and consequently $0\in \sigma_c(T)$. Thus $d(0,\sigma(T))= 0 =\gamma(T)$.\\
Now assume that $\gamma(T)>0$. As a consequence of Theorem \ref{Theorem2}, $\r(T)$ is closed. Note
that $0\notin\sigma(T)$, otherwise Theorem \ref{theorem8} and Theorem \ref{Theorem5} implies that $0\in\sigma_p(T)$.
But this is not true, as $T$ is injective. Thus $0\notin\sigma(T)$ and $T^{-1}$ is bounded totally paranormal
operator, by Theorem \ref{theorem4}. Hence $T^{-1}$ is normaloid and \cite[Proposition 2.12]{KNR}
implies that
\begin{eqnarray*}
  \gamma(T) &=&\frac{1}{\norm{T^{-1}}}=\frac{1}{r(T^{-1})}=\frac{1}{\sup\set{|\mu|:\mu\in\sigma(T^{-1})}} \\
   &=&\inf\set{|\nu|:\nu\in\sigma(T)} =d(0,\sigma(T)).
\end{eqnarray*}
This completes the proof.
\end{proof}
%=================================================================================
As a consequence of Theorem \ref{theorem8} we have the following result.
\begin{corollary}\label{Corollary-A3} If $T\in\lh$ is a densely defined closed totally paranormal operator and
$\n(T)=\n(T^*)$, then $\gamma(T)=d(T):=\inf\set{|\mu|:\mu\in \sigma(T)\setminus\set{0}}$.
\end{corollary}
\begin{proof} Consider the operator $S_0=T|_{\n(T)^{\bot}}:\D(T)\cap\n(T)^{\bot}\to \n(T)^{\bot}$.
 By Theorem \ref{theorem4} and Theorem \ref{theorem8}, $S_0$ is paranormal and
$\gamma(T)=m(S_0)=d(0,\sigma(S_0))=d(T)$.
This proves the result
\end{proof}
%=============================================================
The following example shows the following facts:
\begin{itemize}
  \item  Theorem \ref{theorem8} does not hold if $\n(T)\neq\n(T^*)$.
  \item  It is well known that the residual spectrum of a closed densely defined normal
operator is empty. But this is not true in the case of totally paranormal operators.
\end{itemize}
\begin{example} Let $T:\ell^2(\N)\to\ell^2(\N)$ be defined by
\begin{equation*}
  T(x_1,x_2,\cdots)=(0,x_1,2x_2,3x_3,\cdots),
\end{equation*}
where $\D(T)=\set{(x_1,x_2,\cdots)\in\ell^2(\N):\sum_{j=1}^{\infty}\norm{jx_j}^2<\infty}$.\\
As $c_0$, the space of all complex sequences consisting of at most finitely many
non zero terms is a subset of $\D(T)$ and is dense in $\ell^2(\N)$, we can conclude that $T$ is
densely defined. Hence $T^*$ is well defined. Note that $T$ is a closed operator. We can
show that
\begin{equation*}
  T^*(x_1,x_2,\cdots)=(x_2,2x_3,3x_4,\cdots)
\end{equation*}
with $D(T^*)=\set{(x_n)\in\ell^2(\N):\sum_{j=2}^{\infty}\norm{(j-1)x_j}^2<\infty}$.\\
For any $x=(x_n) \in \D((T-\lambda I)^2)$ and $\lambda\in\c$, we have
\begin{eqnarray*}
  \norm{(T-\lambda I)x}^2 &=&\sum_{j=1}^{\infty}\norm{(j-\lambda)x_j}^2\leq\sum_{j=1}^{\infty}(j+1-\lambda)(j-\lambda)\norm{x_j}^2\\
  &\leq&\bra{\sum_{j=1}^{\infty}\bra{(j+1-\lambda)(j-\lambda)\norm{x_j}}^2}^{\frac{1}{2}}\bra{\sum_{j=1}^{\infty}\norm{x_j}^2} ^{\frac{1}{2}}\\
  &\leq& \norm{(T-\lambda I)^2x}\norm{x}.
\end{eqnarray*}
Hence $T$ is totally paranormal.\\
Since $\norm{Tx}\geq \norm{x}$ for all $x\in \D(T)$ and $\norm{Te_1}=\norm{e_1}$, we get $m(T)=1$. Also
it can be easily verified that $T$ is injective, $\r(T)= \ell^2(\N)\setminus span\set{e_1}$ is closed but
$\r(T)\neq \h$, so $0\in\sigma(T)$. Hence $d(0,\sigma(T))=0\neq 1= m(T)$.\\
\indent Now we will show that $\sigma(T)=\c$. To prove this, we show that $T-\mu I$ is injective
and $\n(T-\mu i)^*\neq \{0\}$, for all $\mu\in\c$.\\
\indent Let $\mu\in\c\setminus\{0\}$ and $(T-\mu I))x=0$ for some $x=(x_n)\in\D(T)$. Then
\begin{equation*}
  \bra{-\mu x_1,x_1-\mu x_2,2x_2-\mu x_3,\cdots}=0.
\end{equation*}
Equating component-wise we get $x=0$. This implies that $T-\mu I$ injective.
Let $y=(y_n) \in \D(T^*)$ be such that $(T-\mu I)^* y =0$. That is,
\begin{equation*}
  \bra{y_2-\bar{\mu}y_1,2y_3-\bar{\mu}y_2,y_4-\bar{\mu}y_3,\cdots}=0.
\end{equation*}
From this we get
\begin{equation}\label{Eq.A-1}
  y=\bra{1,\bar{\mu},\frac{(\bar{\mu})^2}{2!},\frac{(\bar{\mu})^3}{3!},\cdots}y_1.
\end{equation}
If $\mu=0$, then $\n(T^*)= span\{e1\}$. If $\mu\neq 0$, then we will show that $y$ obtained in
Equation (\ref{Eq.A-1}) belongs to $\n(T-\mu I)^*$. Consider $z_n=\frac{\bar{\mu}^{2n}}{(n!)^2}$. Then
\begin{equation*}
  \abs{\frac{z_{n+1}}{z_n}}=\frac{|\mu|^2}{n+1}\to 0\,\,\mbox{as $n\to\infty$}.
\end{equation*}
By the ratio test we conclude that $\sum_{j=1}^{\infty}z_n$ is absolutely convergent, that is
$\sum_{n=1}^{\infty}\bra{\frac{|\mu|^n}{n!}}^2<\infty$. Thus $y\in\ell^2(\N)$. On the similar lines we can show that
$\sum_{n=1}^{\infty}\bra{\frac{|\mu|^n}{(n-1)!}}^2<\infty$. Hence $\n(T-\mu I))^*\neq \{0\}$.\\
For every $\mu\in\c$, $\n(T-\mu I)=\{0\}$ and $\r(T-\mu I)=\bra{\n(T-\mu I)^*}^{\bot}\neq \ell^2(\N)$.
Hence we conclude that $\mu\in\sigma_r(T)$, and $\sigma(T)=\c$.
We also have $\gamma(T)=1\neq 0=d(T)$. From this we conclude that Corollary  \ref{Corollary-A3}
is also not true if the condition, $\n(T)=\n(T^*)$ is dropped.
\end{example}
%==============================================================================================
\begin{theorem}\label{Theorem-Ghad} Suppose $T\in\lh$ is a densely defined closed totally paranormal operator, $\n(T)=\n(T^*)$ and 0 is an isolated spectral
value of $T$ . Then $0 \in\sigma_p(T)$.
\end{theorem}
\begin{proof} Since 0 is an isolated spectral value of $T$ , $d(T )>0$. Hence by Corollary \ref{Corollary-A3},
$\gamma(T)>)$ so that by Theorem \ref{Theorem7}, $\r(T)$ is closed. If $0\notin\sigma_p(T)$, then $\n(T)=\{0\}$
so that we also have $\r(T)= \overline{\r(T)}=\n(T)^{\bot}=\h$, making $T$ bijective and hence
$0\notin \sigma(T)$, a contradiction. Hence $0\in\sigma_p(T)$.
\end{proof}
%======================================================================
\begin{example}  The converse of Theorem \ref{Theorem-Ghad} need not be true. To see this, consider
$T:\ell^2(\N)\to \ell^2(\N)$ defined by
\begin{equation*}
  T(x_1,x_2,x_3,\cdots)=\bra{0,2x_2,\frac{1}{3}x_3,4x_4,\frac{1}{5}x_5,\cdots},
\end{equation*}
where $\D(T)=\set{x\in\ell^2(\N):\bra{0,2x_2,\frac{1}{3}x_3,4x_4,\frac{1}{5}x_5,\cdots}\in\ell^2(\N)}$.
 Here $T$ is a densely defined closed totally paranormal operator.
Since $T$ is not one to one, $0\in \sigma_p(T)$ but it is not an isolated point of the spectrum
$\sigma(T)=\set{0,2,\frac{1}{3},4,\frac{1}{5},\cdots}$.
\end{example}
%====================================================================
\begin{theorem} Suppose $T\in\lh$ is a densely defined closed totally paranormal operator, $\n(T)=\n(T^*)$ . Then $\r(T)$ is closed if and only if 0 is not an accumulation point of $\sigma(T)$.
\end{theorem}
\begin{proof}  By Theorem \ref{Theorem2}, $\r(T)$ is closed if and only if $\gamma(T)>0$ and by Corollary \ref{Corollary-A3},
$\gamma(T)=d(T)$. Hence, $\r(T)$ is closed if and only if $d(T)> 0$ if and only if 0 is not an
accumulation point of $\sigma(T)$.
\end{proof}
%================================================================================
\section{Weyl's theorem for totally paranormal operators}
%====================================================================
In this section, we demonstrate the fulfillment of Weyl's theorem by a densely defined closed operator $T$ that is totally paranormal. Additionally, we establish the self-adjointness of the Riesz projection $E_{\mu}$ corresponding to any non-zero isolated spectral value $\mu$ of $T$.\\
For a Hilbert space $\h$ decomposed as $\h=\h_1\oplus\h_2$, where $T\in\lh$ is a closed operator, we ascertain the block matrix representation of $T$.
\begin{equation}\label{Weyl-eq.1}
  T=\begin{bmatrix}
      T_{11} & T_{12} \\
      T_{21} & T_{22} \\
    \end{bmatrix}
\end{equation}
where $T_{ij}: \D(T)\cap \h_{j}\to \h_i$ is defined by $T_{ij}=P\h_iT P\h_j|_{\D(T)\cap \h_j}$ for $i,j=1,2$.
Here $P_{\h_i}$ is an orthogonal projection onto $\h_i$.\\
For $(x_1,x_2)\in(\h_1\cap \D(T)\oplus(\h_2\cap \D(T))$,
\begin{equation*}
  T(x_1,x_2)=(T_{11}x_1 + T_{12}x_2, T_{21}x_1 + T_{22}x_2).
\end{equation*}
\begin{remark}\label{Remark-Weyl} Let $T$ be as defined in Equation (\ref{Weyl-eq.1}). If $\h_1=\n(T)\neq\{0\}$ and $\h_2=\n(T)^{\bot}$,
then
\begin{equation}\label{Weyl-eq.2}
  T=\begin{bmatrix}
      0 & T_{12} \\
      0 & T_{22} \\
    \end{bmatrix}.
\end{equation}
\begin{itemize}
  \item If $T$ is densely defined closed operator then by Lemma \ref{Lemma1}, $T_{22}\in\L(\n(T)^{\bot})$ is
also densely defined closed operator.
  \item   It can be easily checked that $\r(T_{22})=\r(T)\cap\n(T)^{\bot}$. If $\r(T)$ is closed, then
$\r(T_{22})$ is closed in $\n(T)^{\bot}$.
\end{itemize}
\end{remark}
\indent We say a closed operator $T\in\lh$ satisfy the Weyl's theorem if the Weyl's
spectrum, $\sigma_w(T)$ consists of all spectral values of $T$ except the isolated eigenvalues of
finite multiplicity. That is, $\sigma(T)\setminus\sigma_w(T)=\pi_{00}(T)$.\\
\indent In the work presented by Coburn \cite{Coburn}, it was demonstrated that Weyl's theorem holds true for any bounded hyponormal and Toeplitz operator. Schmoeger \cite{Schmo} later expanded this result to encompass bounded totally paranormal operators. In the current context, we aim to establish the validity of Weyl's theorem for unbounded totally paranormal operators.
%===================================================================================
\begin{theorem}\label{Weyl} Let $T\in\lh$ be a densely defined closed totally paranormal operator. Then
Weyl's theorem holds for $T$, that is, $\sigma(T)\setminus\sigma_w(T)=\pi_{00}(T)$.
\end{theorem}
\begin{proof} Let $\mu\in \sigma(T)\setminus\sigma_w(T)$.  So, we have $\dim(\n(T-\mu I))=\dim(\n(T-\mu I)^{*})<\infty$
and $\r(T-\mu I)$ is closed.\\
On $\h=\n(T-\mu I)\oplus \n(T-\mu I)^{\bot}$, $T-\mu I$ can be decomposed as
\begin{equation*}
  T-\mu I=\begin{bmatrix}
      0 & T_{12} \\
      0 & T_{22}-\mu I|_{\n(T-\mu I)^{\bot}} \\
    \end{bmatrix},
\end{equation*}
where $T_{22}= P|_{\n(T-\mu I)^{\bot}}T|_{\n(T-\mu I)^{\bot}}$. By Remark \ref{Remark-Weyl}, $T_{22}-\mu I_{\n(T-\mu I)^{\bot}}$ is a
densely defined closed operator with domain $\D(T-\mu I)\cap \n(T-\mu I)^{\bot}$ and
$\r(T_{22}-\mu I_{\n(T-\mu I)^{\bot}})$ is closed.\\
As $\n(T-\mu I)$ is finite dimensional, this implies $T_{12}$ is finite rank operator and by
Remark \ref{remark1-intro}, $\ind(T-\mu I)=\ind\bra{T_{22}-\mu I|_{\n(T-\mu I)^{\bot}}}=0$.\\
Since $\n(T_{22}-\mu I|_{\n(T-\mu I)^{\bot}})=\{0\}$ and $\ind\bra{T_{22}-\mu I|_{\n(T-\mu I)^{\bot}}}= 0$, we get
$\n(T_{22}-\mu I|_{\n(T-\mu I)^{\bot}})^{*}=\{0\}$ and consequently $\overline{\r(T_{22}-\mu I|_{\n(T-\mu I)}^{\bot})}= \n(T-\mu I)^{\bot}$.
Thus $T_{22}-\mu I|_{\n(T-\mu I)}^{\bot}$ has bounded inverse and hence $\mu\notin \sigma(T_{22})$. As $\sigma(T)\subseteq \{\mu\}\cup \sigma(T_{22})$, this implies that $\mu$ is an isolated point of $\sigma(T)$. Hence $\mu\in\pi_{00}(T)$.
Conversely, let $\mu\in\pi_{00}(T)$. Now consider the Riesz projection $E_{\mu}$ with respect to
$\mu$. By Theorem \ref{Theorem1} and Theorem \ref{Theorem6}, $\mu\notin \sigma(T|_{\n(E_{\mu})})$  and
$\r(T-\mu I)=\r\bra{\bra{T-\mu I}|_{\n(E_{\mu})}}=\n(E_{\mu})$.\\
Since $\mu\notin \sigma(T|_{\n(E_{\mu})})$, we have that $\r((T-\mu I)|_{\n(E_{\mu})})= \n(E_{\mu})$. Hence $\r(T-\mu I)$
is closed. Also $((T-\mu I)|_{\n(E_{\mu})})^{-1}\in \b\bra{\n(E_{\mu})}$. Thus we get
\begin{eqnarray*}
  \dim\n(T-\mu I)^* &=& \dim\bra{\r(T-\mu I)^{\bot}}=\dim\bra{\n(E_{\mu})^{\bot}} \\
   &=& \dim\bra{\r(E_{\mu})}=\dim\bra{\n(T-\mu I)}.
\end{eqnarray*}
Note that $\dim\bra{\n(E_{\mu})^{\bot}}=\dim\bra{\r(E_{\mu})}$ but the spaces, $\n(E_{\mu})^{\bot}$ and $\r(E_{\mu})$ need not
be the same. Hence $T-\mu I$ is Fredholm operator of index zero. This proves our result.
\end{proof}
%=======================================================================================
\begin{theorem}\label{Theorem-self} Let $T\in\lh$ be a densely defined closed totally paranormal operator and $\mu$
be a non-zero isolated point of $\sigma(T)$. Then the Riesz projection $E_{\mu}$ with respect to $\mu$
satisfy
\begin{equation*}
  \r(E_{\mu})=\n(T-\mu I)=\n(T-\mu I)^*.
\end{equation*}
Moreover, $E_{\mu}$ is self-adjoint.
\end{theorem}
\begin{proof} Let $\mu$ be a non-zero isolated point of $\sigma(T)$. By Theorem \ref{Theorem1} and Theorem \ref{Theorem6},
$\mu\notin \sigma\bra{T|_{\n(E_{\mu})}}$ and $\r(T-\mu I)=\n(E_{\mu})$.
That means $(T-\mu I)|_{\n(E_{\mu})}:\n(E_{\mu})\cap \D(T)\to \n(E_{\mu})=\r(T-\mu I)$ is a bijection.
Also $(T-\mu I)|_{\n(T-\mu I)^{\bot}}\cap\D(T):\n(T-\mu I)\cap\D(T)\to \r(T-\mu I)$ is a bijection,
 we have $\n(E_{\mu})\cap\D(T)\subseteq \n(T-\mu I)^{\bot}\cap\D(T)$.\\
Now we claim that $\n(E_{\mu})\cap\D(T)=\n(T-\mu I)^{\bot}\cap \D(T)$. Let $x\in \n(T-\mu I)^{\bot}\cap \D(T)$ and
$E_{\mu}x=p+q$, where $p\in \n(T-\mu I)$, $q\in \n(T-\mu I)^{\bot}$.
Operating $E_{\mu}$ on both sides, we get
\begin{equation*}
  p+q=E_{\mu}x=p+E_{\mu}q.
\end{equation*}
This implies $E_{\mu}q=q\in\r(E_{\mu})\cap \n(T-\mu I)^{\bot}=\{0\}$, by Theorem \ref{Theorem5}. From this we
conclude that $E_{\mu}x=p=E_{\mu}p$, that is,  $x-p \in \n(E_{\mu})\cap\D(T)\subseteq \n(T-\mu I)^{\bot}\cap\D(T)$.
As $x\in \n(T-\mu I)^{\bot}$, we get $p \in \n(T-\mu I)\cap \n(T-\mu I)^{\bot}=\{0\}$. Consequently
$E_{\mu}x=0$. So $\n(T-\mu I)^{\bot}\cap\D(T)\subseteq \n(E_{\mu})\cap\D(T)$. Hence $\n(T-\mu I)^{\bot}\cap\D(T)=\n(E_{\mu})\cap\D(T)$.
By Lemma \ref{Lemma1} and Theorem \ref{Theorem6}, we get
\begin{eqnarray*}
  \n(T-\mu I)^{\bot} &=&\overline{\n(T-\mu I)^{\bot}\cap\D(T)}=\overline{\n(E_{\mu})\cap\D(T)} \\
   &=& \overline{\r(T-\mu I)\cap\D(T)}=\overline{(\n(T-\mu I)^{*})^{\bot}\cap\D(T)}\subseteq (\n(T-\mu I)^{*})^{\bot}.
\end{eqnarray*}
Hence $\n(T-\mu I)^{*}\subseteq \n(T-\mu I)$. By Theorem \ref{Theorem6},
$\n(E_{\mu})^{\bot}=\r(T-\mu I)^{\bot}=\n((T-\mu I)^*)\subseteq \n(T-\mu I)=\r(E_{\mu})$.
Hence $\n(E_{\mu})^{\bot}=\r(T-\mu I)^{\bot}=\n((T-\mu I)^*)\subseteq \n(T-\mu I)=\r(E_{\mu})$.
 Hence  $\n(E_{\mu})^{\bot}\subseteq \r(E_{\mu})$.
If $x\in \r(E_{\mu})$, then $x=u+v$ where $u\in \n(E_{\mu})$ and $v\in\n(E_{\mu})^{\bot}$. As $\n(E_{\mu})^{\bot}\subseteq \r(E_{\mu})$,
 we get $u=x-v\in \n(E_{\mu})\cap\r(E_{\mu})=\{0\}$. Thus we get $\n(E_{\mu})^{\bot}=\r(E_{\mu})$,
which is equivalent to say that $\n(T-\mu I)=\n(T-\mu I)^*$.
As $\n(E_{\mu})^{\bot}=\r(E_{\mu})$, we have that $E_{\mu}$ is an orthogonal projection. Hence $E_{\mu}$ is
self-adjoint.
\end{proof}
%===================================================================================
From the proof of Theorem \ref{Theorem-self}, we have
\begin{corollary}\label{Mohd19} Let $T\in\lh$ be a densely defined closed paranormal operator.
  Then
  \begin{equation*}
  \n(T-\mu I)^*\subset \n(T-\mu I)
  \end{equation*}
  for all $\mu\in\c$.
\end{corollary}
%=====================================================================================
\indent Using the Birkhoff-James orthogonality concept, we demonstrate that when dealing with a paranormal operator, the eigenspaces associated with distinct isolated eigenvalues are completely independent of each other. To explain, consider a subspace $M$ within a Banach space $X$. We say $M$ is Birkhoff–James orthogonal to another subspace $N$ in $X$ if, $\norm{m}\leq \norm{m+n}$ for all $m\in M$ and $n\in N$.  This definition aligns with the Birkhoff–James orthogonality. When $X$ is a Hilbert space, this notion coincides with the standard idea of orthogonality.
%==========================================================================================
\begin{proposition}\label{perpendicular} Let $T\in\lh$ be a densely defined closed totally paranormal operator. If $\mu_1$
and $\mu_2$ are two non zero distinct isolated points of $\sigma(T)$, then $\n(T-\mu_1 I)$ is orthogonal
to $\n(T-\mu_2 I)$.
\end{proposition}
\begin{proof} Without loss of generality, assume that $|\mu_1|<|\mu_2|$. For any $x\in\n(T-\mu_1 I)$
and $y\in \n(T-\mu_2 I)$, consider the set $M=span\{x, y\}$. As $M$ is invariant subspace
for $T$ , it follows that $T|_{M}$ is totally paranormal operator and $\norm{T|_{M}}=|\mu_2|$. We have the
following.
\begin{eqnarray*}
 \norm{\frac{\mu_1^n}{\mu_2^n}x+y} &=& \frac{1}{\abs{\mu_2^n}}\norm{\mu_1^n x+\mu_2^n y} \\
   &\leq&\frac{\norm{T|_{M}}^n}{\abs{\mu_2^n}}\norm{x+y}=\norm{x+y}
\end{eqnarray*}
Taking the limit $n\to\infty$, we get $\norm{y}\leq \norm{x+y} $, for every $x\in \n(T-\mu_1I)$ and
$y\in\n(T-\mu_2 I)$. Hence $\n(T-\mu_2 I)$ is orthogonal to $\n(T-\mu_1I)$.\\
Next, if $|\mu_1|=|\mu_2|$, then for every $n\in \N$
\begin{eqnarray*}
  \norm{\bra{\frac{\mu_1+\mu_2}{2\mu_2}}^nx+y} &=& \norm{\frac{(\mu_1+\mu_2)^nx+(\mu_1+\mu_2)^ny}{(2\mu_2)^n}} \\
   &\leq& \frac{1}{(2|\mu_2|)^n}\sum_{j=0}^{n}\binom{n}{j}|\mu_2|^{j}\norm{\mu_1^{n-j}x+\mu_2^{n-j}y}\\
   &=&  \frac{1}{(2|\mu_2|)^n}\sum_{j=0}^{n}\binom{n}{j}|\mu_2|^{j}\norm{\bra{T|_{M}}^{n-j}(x+y)}\\
   &\leq&  \frac{1}{(2|\mu_2|)^n}\sum_{j=0}^{n}\binom{n}{j}|\mu_2|^{n}\norm{x+y}\\
   &=&\norm{x+y}.
\end{eqnarray*}
As $\mu_1\neq \mu_2$, we have $\abs{\frac{\mu_1+\mu_2}{2\mu_2}}<1$. Now as $n\to\infty$ in the above inequality we get
that $\norm{y}\leq \norm{x+y}$. This proves the result
\end{proof}
%==========================================================================
%=============================================================================
\begin{proposition}
 Let $T\in\lh$ be a densely defined closed totally paranormal operator and  $T^{2}$ be a compact operator. Then $T$
  is also compact and normal.
\end{proposition}
\begin{proof}
  Assume that $T$ is a totally paranormal operator . Hence,
  \begin{equation}\label{M2}
    \norm{Tx}^2\leq \norm{T^{2}x}\norm{x}\,\,\,\mbox{for every $x\in\D(T^2)$}.
  \end{equation}
   Let $\{x_m\}\in\h$ be weakly convergent sequence with limit $0$ in $\D(T)$. From the compactness
  of $T^2$ and the relation (\ref{M2}) we get the following relation:
  $$ \norm{Tx_m}^{2}\rightarrow 0,\,\,m\rightarrow \infty.$$
  From the last relation it follows that $T$ is compact.
  Since $T$ is compact $\sigma(T)$ is finite set or countable infinite with 0 as the unique
limit point of it. Let $\sigma(T) \setminus \{0\} = \{\lambda_n\}$ with
$$|\lambda_1|\geq |\lambda_2|\geq\cdots \geq |\lambda_n|\geq |\lambda_{n+1}|\geq \cdots\geq 0,\,\,\, \mbox{and}\,\,
\lambda_n\rightarrow 0\,\,(n\rightarrow \infty).$$
By the compactness of $T$ or isoloidness of $T$, $\lambda_n \in \sigma_p(T)$ and $dim \ker(T-\lambda_n) <\infty$
for all $n$. Since $\ker(T-\lambda_n)\subset \ker(T-\lambda_n)^*$, $\M:=\bigoplus_{n=1}^{\infty}\ker(T-\lambda_n)$
reduces $T$, and $T$ is of the form
$$T=\bra{\bigoplus_{n=1}^{\infty}\lambda_n}\oplus T'\,\,\mbox{on}\,\, \h=\M\oplus \M^{\bot}.$$
By the construction, $T'$ is totally paranormal and $\sigma(T')=\{0\}$ hence $T'=0$. This shows
that
$$T=\bra{\bigoplus_{n=1}^{\infty}\lambda_n}\oplus 0$$
and it is normal.
\end{proof}
%===============================================================================
%=================================================================
\begin{theorem} Let $T\in\lh$ be a densely defined closed totally paranormal operator with $\sigma_w(T)=\{0\}$. Then
 $T$ is a compact normal operator.
\end{theorem}
\begin{proof} By Theorem \ref{Weyl}, $T$ satisfy Weyl's theorem and this implies that
each element in $\sigma(T)\setminus\sigma_w(T)=\sigma(T)\setminus\{0\}$ is an eigenvalue of $T$ with finite multiplicity,
and is isolated in $\sigma(T)$. Hence $\sigma(T)\setminus\{0\}$ is a finite set or a countable set with 0 as its only accumulation point.
Put $\sigma(T)\setminus\{\lambda_n\}$, where $\lambda_n\neq \lambda_m$ whenever $n\neq m$ and $\{|\lambda_n|\}$ is a non-increasing sequence.
Since $T$ is normaloid, we have $|\lambda_1|=\norm{T}.$  By Corollary \ref{Mohd19}, we have $(T-\lambda_1 I)x=0$ implies $(T-\lambda_1 I)^*x=0$.
Hence $\n(T-\lambda_1 I)$ is a reducing subspace of $T$. Let $E_1$ be the orthogonal projection onto $\n(T-\lambda_1 I).$
Then $T=\lambda_1 I\oplus T_1$ on $\h=\r(E_1)\oplus \r(I-E_1)$. Since $T_1$ is totally paranormal by Theorem \ref{theorem4} (i) and $\sigma_p(T)=\sigma_p(T_1)\cup\{\lambda_1\}$,
we have $\lambda_2\in \sigma_p(T_1)$. By the same argument as above, $\n(T-\lambda_2 I)=\n(T_1-\lambda_2I)$ is a finite dimensional reducing subspace
of $T$ which is included in $\r(I-E_1)$. Put $E_2$ be the othogonal projection onto $\n(T-\lambda_2 I)$. Then
$T=\lambda_1E_1\oplus \lambda_2 E_2\oplus T_2$ on $\h=\r(E_1)\oplus \r(E_2)\oplus \r(I-E_1-E_2)$. By repeating above argument,
 each $\n(T-\lambda_n I)$ is a reducing subspace of $T$ and $\norm{T-\displaystyle{\bigoplus_{k=1}^{n}\lambda_kE_k}}=\norm{T_n}=|\lambda_{n+1}|\to 0$
 as $n\to \infty$. Here $E_k$ is the orthogonal projection onto $\n(T-\lambda_k)$ and $T=(\displaystyle{\bigoplus_{k=1}^{n}\lambda_kE_k})\oplus T_n$ on
 $\h=\displaystyle{\bigoplus_{k=1}^{n}\r(E_k)}\oplus (1-\displaystyle{\sum_{k=1}^{n}\r(E_k)}.$ Hence $T=\displaystyle{\bigoplus_{k=1}^{\infty}\lambda_kE_k}$
  is compact and normal because each $E_k$ is a finite rank orthogonal projection which satisfies $E_kE_t=0$ whenever $k\neq t$ by
  Proposition \ref{perpendicular} and $\lambda_n\to 0$ as $n\to \infty$.
\end{proof}
%==============================================================================
%=====================================================================
\section{Conclusion and future Work}
In conclusion, this article has delved into a comprehensive analysis of spectral properties related to totally paranormal closed operators within the context of Hilbert spaces. The exploration extended beyond the conventional bounds of boundedness, incorporating closed symmetric operators into the discussion.

The initial focus was on establishing the non-emptiness of the spectrum for such operators, accompanied by a characterization of closed-range operators based on the spectrum. Building on these foundational results, Weyl's theorem was proven for densely defined closed totally paranormal operators. Specifically, it was demonstrated that the difference between the spectrum $\sigma(T)$ and the Weyl spectrum $\sigma_w(T)$ is precisely the set of isolated eigenvalues with finite multiplicities, denoted as $\pi_{00}(T)$.

The final section of the article explored the self-adjointness of the Riesz projection $E_{\mu}$ corresponding to any non-zero isolated spectral value $\mu$ of the operator $T$. The relationships $\mathrm{ran}(E_{\mu}) = \n(T-\mu I) = \n(T-\mu I)^*$ were established for this Riesz projection. Furthermore, it was shown that if a closed totally paranormal operator $T$ has a Weyl spectrum $\sigma_w(T) = {0}$, then $T$ qualifies as a compact normal operator.

In terms of future work, potential avenues include exploring applications of these spectral properties in specific mathematical or physical contexts. Additionally, investigating the implications of these results on related areas of operator theory or functional analysis could provide valuable insights. Further developments in the understanding of totally paranormal operators and their spectral characteristics may contribute to advancements in various mathematical disciplines.
%========================================================================== ====

\bibliographystyle{unsrtnat}
\bibliography{references}  %%% Uncomment this line and comment out the ``thebibliography'' section below to use the external .bib file (using bibtex) .




% 	\bibitem{keshet2016prediction}
% 	Keshet, Renato, Alina Maor, and George Kour.
% 	\newblock Prediction-Based, Prioritized Market-Share Insight Extraction.
% 	\newblock In {\em Advanced Data Mining and Applications (ADMA), 2016 12th International
%                       Conference of}, pages 81--94,2016.

% \end{thebibliography}

\end{document}